\documentclass[twoside,11pt]{article}

\usepackage{psfrag}
\usepackage{graphicx}
\usepackage{amssymb}
\usepackage{latexsym}
\usepackage{amsmath}

\newtheorem{theorem}{Theorem}
\newtheorem{proposition}{Proposition}

\newtheorem{lemma}{Lemma}
\newtheorem{corollary}{Corollary}

\newenvironment{proof}[1][Proof]{\noindent\textbf{#1} }{\ \rule{0.5em}{0.5em}}%preliminary
\date{}

\begin{document}

\title{Hilbert metrics and Minkowski norms}
\maketitle

\begin{center}
{\large Thomas Foertsch$^\dagger$ \hspace{1cm} Anders
Karlsson$^\ddagger$} \footnote{$\dagger$ Department of Mathematics, University of Michigan,
525 East University, Ann Arbor, USA, E-mail: foertsch@math.unizh.ch}
\footnote{$\ddagger$ Department of Mathematics, Royal Institute
of Technology, 100 44 Stockholm, Sweden,
E-mail: akarl@math.kth.se}
 \\
\end{center}

\begin{abstract}
It is shown that the Hilbert geometry $(D,h_D)$
associated to a bounded convex domain $D\subset \mathbb{E}^n$ is
isometric to a normed vector space $(V,||\cdot ||)$ if and only if
$D$ is an open $n$-simplex. One further result on the asymptotic
geometry of Hilbert's metric is obtained with corollaries
for the behavior of geodesics. Finally we prove that every geodesic ray in 
a Hilbert geometry converges to a point of the boundary.
\end{abstract}

\vspace{0.5cm}

%%%%%%%%%%%%%%%%%%%%%%%%%%%%%%%%%%%%%%%%%%%%%%%%%%%%%%%%%%%%%%%%%%%%%%%%%%%%%%%%%%%%%%%%%%%%%
%%%%%%%%%%%%%%%%%%%%%%%%%%%%%%%%%%%%%%%%%%%%%%%%%%%%%%%%%%%%%%%%%%%%%%%%%%%%%%%%%%%%%%%%%%%%%

\section{Introduction}
Busemann wrote on page 105 in \cite{b} that ``Plane Minkowskian
geometry arises from the Euclidean through replacing the ellipse
as unit circle by a convex curve. In a somewhat similar way a
geometry discovered by Hilbert arises from Klein's Model of
hyperbolic geometry through replacing the ellipse as absolute
locus by a convex curve." In this note we treat the question of
when a Hilbert geometry is isometric to a Minkowski space (here
meaning a normed finite dimensional real vector space).

We recall the definition of Hilbert's metric. Let $\mathbb{E}^n$
denote the $n$-dimensional Euclidean space. For the Euclidean
distance of $x,y\in \mathbb{E}^n$ we write $|xy|$, for the
straight line segment between $x$ and $y$ we write $[x,y]$ and
$L(x,y)$ denotes the whole straight line through $x$ and $y$. \\
Given a bounded convex domain $D\subset \mathbb{E}^n$ with boundary $\partial D\subset \mathbb{E}^n$, the
Hilbert metric $h_D:D\times D\longrightarrow \mathbb{R}_0^+$ is defined as follows: For
$x,y\in D$ one defines
\begin{displaymath}
h_D (x,y) \; := \;
\left\{
\begin{array}{cl}
\log \; \frac{|y\bar{x}|\cdot |x\bar{y}|}{|x\bar{x}|\cdot |y\bar{y}|} & \mbox{if} \; x\neq y \\
0 & \mbox{if} \; x=y
\end{array}
\right\} ,
\end{displaymath}
where $\bar{x}\in L(x,y)\cap \partial D$ with
$|\bar{x}x|<|\bar{x}y|$ and $\bar{y}\in L(x,y)\cap \partial D$
with $|\bar{y}x|>|\bar{y}y|$. The expression
$\frac{|y\bar{x}|\cdot |x\bar{y}|}{|x\bar{x}|\cdot |y\bar{y}|}$ is
called the cross ratio of the four collinear ordered points
$\bar{x},x,y,\bar{y}$ and is invariant under projective
transformations. For the basic properties of the distance
$h_D$ see \cite{b} and \cite{dlh}. We shall refer to $(D,h_D)$ as a
\emph{Hilbert geometry}. All vector spaces are over the reals.

The following fact is known, see the pages 22-23 in Nussbaum's
book \cite{nuss} as well as the pages 110-111 and 113 in de
la Harpe's paper \cite{dlh}:
\begin{theorem} \label{lemma-dlh} %\cite{dlh}
Let $D\subset \mathbb{R}^n$ be the interior of the standard
$n$-simplex. Then $(D,h_D)$ is isometric to an
$n$-dimensional normed vector space.
\end{theorem}

%That author gives an explicit isometry and provides the proof in
%dimension $n=2$.
Recall that a \emph{polytope} is the convex hull
of finite number of points and an $n$-\emph{simplex} in
$\mathbb{R}^n$ is the convex hull of $n+1$ points in general
position (hence it has nonempty interior).

In Section \ref{sec-proof} of this paper we show that the 
converse of Theorem \ref{lemma-dlh}
is also true, so that we have:
\begin{theorem} \label{theo-main}
A Hilbert geometry $(D,h_D)$ is isometric to a
normed vector space if and only if $D$ is the interior of a
simplex.
\end{theorem}

%% standard simplex versus simplex?
%%
In particular, a Hilbert geometry is never a Hilbert space. 
%In fact even more is true:
%\begin{theorem} \label{theo-flat}
%Every Hilbert geometry has euclidean rank $1$, or in other words
%it contains no $2$-flat.
%\end{theorem}

In Section \ref{sec-asympt} we obtain an 
estimate of the asymptotic geometry of Hilbert's metric (see Proposition \ref{prop-converse}), which will be used to prove:

\begin{theorem} \label{theo-convergence}
Let $(D,h_D)$ be a Hilbert geometry. Then
\begin{description}
\item[(i)] every geodesic ray in $(D,h_D)$ has to converge to a point in $\partial D$.
\item[(ii)] every complete geodesic in $(D,h_D)$ has precisely two accumulation points on $\partial D$.
\end{description}
\end{theorem}

Theorem \ref{theo-convergence} in the case $D$ is the interior of a simplex
was already proved by V. Metz with entirely different techniques 
(\cite{met}).

%%%%%%%%%%%%%%%%%%%%%%%%%%%%%%%%%%%%%%%%%%%%%%%%%%%%%%%%%%%%%%%%%%%%%%%%%%%%%%%%%%%%%%%%%%%%%%%%%%%%%%%%%%%%%%%%%
%%%%%%%%%%%%%%%%%%%%%%%%%%%%%%%%%%%%%%%%%%%%%%%%%%%%%%%%%%%%%%%%%%%%%%%%%%%%%%%%%%%%%%%%%%%%%%%%%%%%%%%%%%%%%%%%%

\section{The proof of Theorem \ref{theo-main}}
\label{sec-proof}
As is explained both in \cite{b} and \cite{dlh}, the straight line
segments in $D$ are geodesics in $(D,h_D)$, but in
general there can exist geodesics different from such Hilbert
straight line segments. Indeed, two points $x,y\in D$ can be
joined by no other geodesic than their connecting Hilbert straight
line segment, if and only if there do not exist two coplanar but
not collinear straight line segments $l_1,l_2\subset
\partial D$ through $\bar{x}$ and $\bar{y}$ such that $\bar{x}$
and $\bar{y}$ are not boundary points of $l_1$
and $l_2$. \\

%% Straightline geodesics vs geodesics?
%%
The following notation is by now standard:
$$
(x|y)_{p_0}:=\frac{1}{2}\Big[ h_D(x,p_0)+h_D(y,p_0)-h_D(x,y)\Big] ,
$$
where $x$, $y$, and $p_0$ are three points in the metric space. We
recall the following simple but useful fact:
\begin{theorem} \label{theo-anders} (\cite{kn})
Let $D$ be a bounded convex domain. Let
$\{x_n\}$, $\{z_n\}$ be two sequences of points in $D$. Assume
that $x_n\longrightarrow \bar{x}\in \partial D$,
$z_n\longrightarrow \bar{z}\in \partial D$ and $[\bar{x},\bar{z}]
\not\subset \partial D$. Then, for any fixed $p_0$, there is a
constant $K=K(p_0,\bar{x},\bar{z})$ such that
\begin{displaymath}
\limsup_{n\longrightarrow \infty} (x_n|z_n)_{p_0} \; \le \; K.
\end{displaymath}
\end{theorem}

We will need:

\begin{lemma} \label{lemma-main}
Let $D\subset \mathbb{R}^m$ be a bounded convex domain such that
$(D,h_D)$ is isometric to an $m$-dimensional normed
vector space $(V,||\cdot ||)$ and such that there exist
$\bar{x}_i\in \partial D$, $i=1,...,n$ with
$[\bar{x}_i,\bar{x}_j]\not\subset \partial D$ for $i\neq j$. Then
there exist $n$ points $v_1,...,v_n$ on the unit sphere of
$(V,||\cdot ||)$ such that
\begin{equation} \label{eqn-cor-main}
||v_i - v_j|| = 2 \hspace{1cm} \forall i,j=1,...,n, \; i\neq j.
\end{equation}
\end{lemma}

\begin{proof}: Let $\varphi :D\longrightarrow V$ be an isometry with
$\varphi (p_0)=0$. Let further ${\gamma}_i:[0,\infty
)\longrightarrow D$ be the arc length parameterized Hilbert
straight line geodesic connecting $p_0$ to $\bar{x}_i$. In view of
Theorem \ref{theo-anders} we find $\tilde{K}>0$ and $k_0>0$ such
that
\begin{displaymath}
h_D\Big( {\gamma}_i(k),{\gamma}_j(k)\Big) \; \ge \; 2k-2\tilde{K}
\hspace{1cm} \forall i,j=1,...,n, \;\; i\neq j,
\end{displaymath}
for all $k\geq k_0$. Given $N>0$, let
$k_N:=\max\{k_0,2\tilde{K}N\}$. 
For each $i$ let 
$$
v^N_i=\frac {1}{k_N}\varphi(\gamma_i(k_N))
$$
which is a point on the unit sphere in $(V,||\cdot ||)$. Then we have
that
\begin{displaymath}
||v^N_i-v^N_j||=\frac {1}{k_N}
h_D({\gamma}_i(k_{N}),{\gamma}_j(k_{N}))\ge
\frac{1}{k_N}(2k_{N}-2\tilde{K}) \ge  2 - \frac{1}{N}.
\end{displaymath}

%Denote by ${\eta}_j:[0,\infty
%)\longrightarrow V$ the straight line geodesic in $(V,||\cdot ||)$
%such that ${\eta}_j(0)=0$ and
%${\eta}_j(k_{N})=\varphi({\gamma}_{j}(k_{N}))$, $j=1,...,n$. Then for all
%$i,j=1,...,n$, $i\neq j$ we find
%\begin{displaymath}
%||{\eta}_i(1)-{\eta}_j(1)|| =  \frac{d(
%{\gamma}_i(k_{N}),{\gamma}_j(k_{N}))}{k_{N}} \ge
%\frac{2k_{N}-2\tilde{K}}{k_{N}} \ge  2 - \frac{1}{N}.
%\end{displaymath}
Since $V$ is finite dimensional, the unit sphere is compact and
hence we can find a subsequence $N_k \rightarrow\infty$ and
$v_1,...,v_n$ such that $v^{N_k}_i\rightarrow v_i$ for every $i$.
These limit points clearly satisfy (\ref{eqn-cor-main}).
\end{proof}

\begin{proposition}\label{prop-polytope}
Suppose that a Hilbert metric space $(D,h_D)$ is
isometric to a normed space. Then $D$ is the interior of a
polytope.
\end{proposition}
\begin{proof}: Suppose that $\overline D$ is not the convex hull of
a finite number of points. Then one can find an infinite number of
points satisfying the hypothesis of Lemma \ref{lemma-main}.
To see that assume we can only find a finite number of such points.
Take a maximal such set of points of the boundary. Note that none 
of these points belongs to the boundary of two different faces,
because otherwise we could replace this point by two interior
points of these faces, which contradicts the maximality of the
chosen set. Now take the union of the closed faces containing
points in our maximal set (a priori such a face might just be 
a point itself). If this union is not all of the boundary, then we
add a point outside this union to our chosen set, hence again
contradicting maximality. We have thus showed that $\overline D$ 
is the convex hull of a finite number of points.

Therefore if $\overline D$ is not a polytope but isometric to 
a normed vector space, then by Lemma \ref{lemma-main} and a 
diagonal process we can extract an infinite sequence of 
$v_i$ of mutual distance $2$. This clearly contradicts the 
compactness of the unit sphere in $(V,||\cdot ||)$.
\end{proof} \\

With Proposition \ref{prop-polytope} at hand, we are ready to provide the proof of Theorem \ref{theo-main}: \\
\begin{proof}: The ``if part'' follows
from Theorem \ref{lemma-dlh}. For the ``only if part'' we know by
Proposition \ref{prop-polytope} that $\overline{D}$ must be a
polytope. Suppose that $\overline{D}$ is not an $n$-simplex. Then we
find three points $v_1,v_2,e\in \partial D$ such that $v_1$ and
$v_2$ are vertices of $\partial D$, the Hilbert straight line
geodesics $[e,v_1]$ and $[e,v_2]$ lie entirely in $D$ and the
intersection $\partial C$ of the affine plane $\Sigma$ through
$e,v_1,v_2$ is a polytope such that $e$ is not a vertex point of
$\partial C$. We also write $C:=\Sigma \cap D$. Note that $h_C=h_D|_C$.\\
Under the isometry $\varphi : (D,h_D)\longrightarrow (V,||\cdot ||)$ the Hilbert straight line
geodesics $[e,v_1]$ and $[e,v_2]$ have to be mapped to two
straight line geodesics $l_1:=\varphi ([e,v_1])$ and $l_2:=\varphi
([e,v_2])$ in $(V,||\cdot ||)$ because of the uniqueness of these
geodesics. We will get a contradiction to the
assumption that $\partial D$ is not an $n$-simplex by showing that (i) $l_1\parallel l_2$ and (ii) $l_1\not\parallel l_2$: \\
(i) Let $e_1,e_2\not\in \{ v_1,v_2\}$ be the vertices of $\partial C$ with $e\in [e_1,e_2]$ such that
$e_i$ and $v_i$ lie in the same connected  component of $\bar{C}\setminus [v_1,e]$. Let further $\tilde{e}_i\neq e_j$
be the vertex point of $\partial C$ next to $e_i$, $i=1,2$, $i\neq j$, consider a straight line segment $s$ parallel to
$[e_1,e_2]$ with endpoints $\tilde{x}_1$ and $\tilde{x}_2$ on $[e_1,\tilde{e}_1]$ and $[e_2,\tilde{e}_2]$, respectively.
By $x_1$ and $x_2$ we denote the intersection of $s$ with $[e,v_1]$ and $[e,v_2]$, respectively. \\
All we need to show is that there exists $c\in \mathbb{R}^+$ such that (a) for all $y_1\in [e,x_1]$ there exists
$y_2\in [e,x_2]$ satisfying $h_D(y_1,y_2)\le c$ and (b) for all
$y_2\in [e,x_2]$ there exists  $y_1\in [e,x_1]$ satisfying $h_D(y_1,y_2)\le c$, as this obviously implies
$l_1\parallel l_2$. \\
Without loss of generality we only prove (a). In order to do that, we set
$\tilde{x}_1:=s\cap [e_1,v_1]$ and $\tilde{x}_2:=s\cap [e_2,v_2]$. Given $y_1\in [e,x_1]$ we denote by $r$ the straight
line segment through $y_1$ in $D$ parallel to $[e_1,e_2]$ and set $y_2:=r\cap [e,v_2]$. Then
\begin{displaymath}
h_D(y_1,y_2) \; \le \; \log \, \Big[
1+\frac{|x_1x_2|}{|x_1\bar{x}_1|} +
\frac{|x_1x_2|}{|x_2\bar{x}_2|} +
\frac{|x_1x_2|^2}{|x_1\bar{x}_1||x_2\bar{x}_2|}\Big] \; =: \; c,
\end{displaymath}
where $\bar{x}_i:=s\cap [e,v_i]$.

(ii) Let $\tilde{v}_i\neq v_j$ be the vertex of $\partial C$ such
that $[\tilde{v}_i,v_i]\subset
\partial C$, $i,j=1,2$, $i\neq j$. Given $z_i\in C$ we denote by
$\tilde{z}_i\in \partial C$ the unique point on $\partial C$ such
that $z_i\in [\tilde{z}_i,v_j]$, $i,j=1,2$, $i\neq j$. Let $u_i\in
[v_i,e]$ be such that $[u_1,u_2]\parallel [v_1,v_2]$ and
$L(u_1,u_2)\cap \partial C \subset [v_1\tilde{v}_1]\cup
[v_2,\tilde{v}_2]$, $i=1,2$. Let further ${\gamma}_i:(-\infty
,\infty )\longrightarrow D$ be the arc length parameterization of
$[e,v_i]$
with $[{\gamma}_i(-\infty )]=e$ and $[{\gamma}_i(\infty )]=v_i$, $i=1,2$. \\
For $t_0\in \mathbb{R}$ we write $w_1:={\gamma}_1(t_0)$ and choose $s_0\in \mathbb{R}$ such that
\begin{displaymath}
h_D\Big( w_1,{\gamma}_2(s_0)\Big) \; = \; \inf\limits_{s\in \mathbb{R}} h_D\Big( w_1,{\gamma}_2(s)\Big) .
\end{displaymath}
Given $\tilde{c}$ arbitrary large we can choose $z_i\in tr({\gamma}_i)$ such that $z_i\in [v_i,u_i]$,
$\tilde{z}_i\in [v_i,\tilde{v}_i]$ and
$|z_i\tilde{z}_i|<e^{-\frac{\tilde{c}}{2}}|u_1u_2|$, $i=1,2$. \\
Now there exists $t\in \mathbb{R}$ such that 
for $w_1:=[{\gamma}_1(t_0+t),v_2]\cap [v_1,\tilde{v}_1]$ and $w_2:=[{\gamma}_2(s_0+t),v_1]\cap [v_2,\tilde{v}_2]$ it holds
\begin{displaymath}
\max \Big\{ |{\gamma}_i(t_i+t)v_i|,|{\gamma}_i(t_i+t)w_i|\Big\} \; < \; |z_i\tilde{z}_i| \hspace{1cm} i=1,2,
\end{displaymath}
where $t_1:=t_0$ and $t_2:=s_0$. For such a $t$ it is easy to see that
\begin{displaymath}
h_D\Big( {\gamma}_1(t_0+t), {\gamma}_2(s_0+t)\Big) \; \ge \; \tilde{c} ,
\end{displaymath}
which contradicts $l_1\parallel l_2$.
\end{proof}

\begin{figure}[htbp]
\psfrag{v1}{$v_1$} \psfrag{v2}{$v_2$}
\psfrag{tildev1}{$\tilde{v}_1$} \psfrag{tildev2}{$\tilde{v}_2$}
\psfrag{e}{$e$} \psfrag{e1}{$e_1$} \psfrag{e2}{$e_2$}
\psfrag{tildee1}{$\tilde{e}_1$} \psfrag{tildee2}{$\tilde{e}_2$}
\psfrag{x1}{$x_1$} \psfrag{x2}{$x_2$} \psfrag{y1}{$y_1$}
\psfrag{y2}{$y_2$} \psfrag{tildex1}{$\tilde{x}_1$}
\psfrag{tildex2}{$\tilde{x}_2$} \psfrag{barx1}{$\bar{x}_1$}
\psfrag{barx2}{$\bar{x}_2$} \psfrag{u1}{$u_1$} \psfrag{u2}{$u_2$}
\psfrag{z1}{$z_1$} \psfrag{z2}{$z_2$}
\psfrag{tildez1}{$\tilde{z}_1$} \psfrag{tildez2}{$\tilde{z}_2$}
\psfrag{g1(t0)}{${\gamma}_1(t_0)$}
\psfrag{g2(s0)}{${\gamma}_2(s_0)$}
\psfrag{g1(t0+t)}{${\gamma}_1(t_0+t)$}
\psfrag{g2(s0+t)}{${\gamma}_2(s_0+t)$}
\includegraphics[width=0.9\columnwidth]{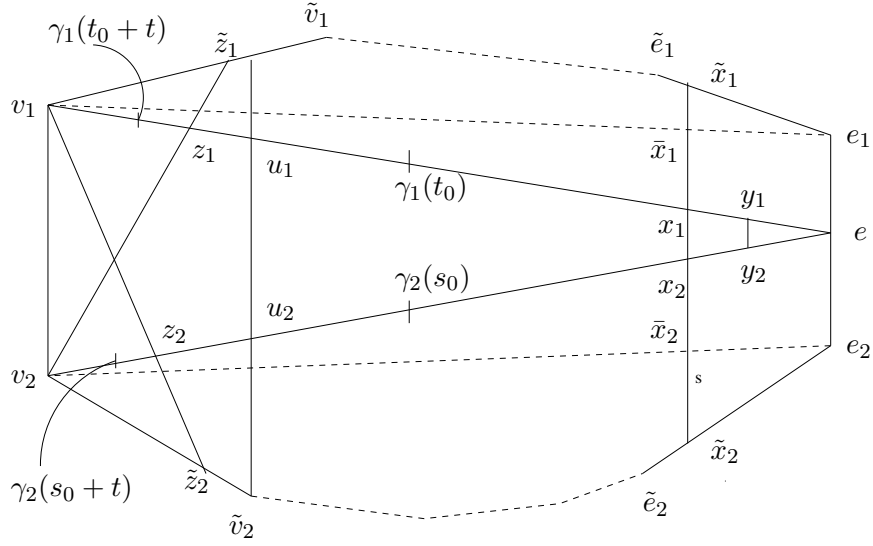}
\caption{The figure illustrates the situation in the proof of
Theorem \ref{theo-main}.}
\end{figure}

\section{Some asymptotic geometry} \label{sec-asympt}

We will need:
\begin{lemma} \label{lemma-angle}
Let $D$ be an open bounded convex domain in $\mathbb{E}^n$ and $p_0\in D$. Then there exists $C=C(D,p_0)>0$ such that
the following holds: Let $z\in \partial D$, $z'\in [p_0,z]\cap D$ and $x,y\in \partial D$ such that $z'\in [x,y]$, then
\begin{equation} \label{eqn-a/b}
|zz'| \; \le \; C|z'x|. 
\end{equation}
\end{lemma}
\begin{proof}:
Without loss of generality assume $x\neq z\neq y$ and let $\delta > 0$ be such that the Euclidean ball $B(p_0,3\delta )$
of radius $3\delta$ around $p_0$ is contained in $D$. Let further $\Sigma$ denote the intersection of $D$ with the affine
plane spaned by $x$, $y$ and $z$, and set $\tilde{B}:=B(p_0,\delta )\cap D$. \\
Let $\gamma$ be the angle between the two straight lines $T_1$ and $T_2$ tangent to $B$ with $T_1\cap T_2=\{ p_0\}$. Since $D$
is bounded, there exists some ${\gamma}_0>0$ such that for all $z\in D$ it holds $\gamma \ge {\gamma}_0$. \\
If $|z'x|\ge \delta$, then inequality \ref{eqn-a/b} holds for $C:=\frac{diam D}{\delta}$. Let now 
$|z'x|< \delta$. Let $\alpha := \angle ([xz],[xz'])$ and $\beta := \angle ([zz'],[zx])$. Then by the sine law we find
\begin{displaymath}
|zz'| \; \le \; \frac{\sin \alpha}{\sin \beta} |z'x|.
\end{displaymath}
Since $|z'x|<\delta$, we also have
\begin{displaymath}
\frac{{\gamma}_0}{2} \; \le \; \beta \; \le \; \pi \; - \; \frac{{\gamma}_0}{2}
\end{displaymath}
and the claim follows.
\end{proof}

%\begin{figure}
%  \includegraphics{interchords.eps}\\
%  \caption{Intersecting chords}\label{Fi:xxx}
%\end{figure}

\begin{figure}[htbp]
\psfrag{x}{$x$} 
\psfrag{y}{$y$}
\psfrag{z}{$z$}
\psfrag{alpha}{$\alpha$}
\psfrag{beta}{$\beta$}
\psfrag{gamma}{$\gamma$}
\psfrag{p0}{$p_0$}
\psfrag{T1}{$T_1$}
\psfrag{T2}{$T2$}
\psfrag{Sigma}{$\Sigma$}
\psfrag{tildeB}{$\tilde{B}$}
\psfrag{z'}{$z'$}
\psfrag{z}{$z$}
\psfrag{B(p0,3delta)}{$\scriptstyle{B(p_0,3\delta )}$}
\psfrag{xbar}{$\bar{x}$}
\psfrag{bn}{$\scriptstyle{b_n}$}
\psfrag{an}{$\scriptstyle{a_n}$}
\psfrag{cn}{$\scriptstyle{c_n}$}
\psfrag{dn}{$\scriptstyle{d_n}$}
\psfrag{en}{$\scriptstyle{e_n}$}
\psfrag{xn}{$\scriptstyle{x_n}$}
\psfrag{yn}{$\scriptstyle{y_n}$}
\psfrag{p0}{$\scriptstyle{p_0}$}
\includegraphics[width=0.9\columnwidth]{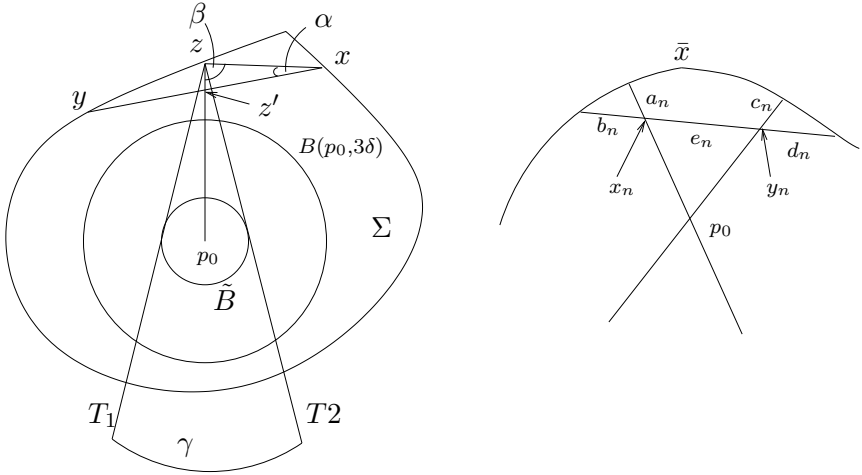}
\caption{The figure on the left hand side illustrates the situation in the proof of
Lemma \ref{lemma-angle}, while the figure on the right hand side explains the notation in the proof of
Proposition \ref{prop-converse}.}
\end{figure}

The following estimate complements Theorem \ref{theo-anders}
above:

\begin{proposition} \label{prop-converse}
Let $\{x_{n}\}$ and $\{y_{n}\}$ be two sequences of points in $D$
both converging to $\bar{x}\in\partial D$. Then for any fixed
$p_0$
$$
(x_{n}|y_{n})_{p_0}\rightarrow\infty.
$$
\end{proposition}
%\begin{figure}
%  \includegraphics{xnynp0.eps}\\
%  \caption{Proposition}\label{Fi:xxx}
%\end{figure}

\begin{proof}:
See Figure 2 for the notation. We can assume that
$d(x_n,y_n)\rightarrow\infty$ (because for any subsequence for
which $x_n$ and $y_n$ stay bounded the conclusion for that
subsequence is immediate). Therefore at least one of $b_n$ or
$d_n$ tend to $0$.

By compactness of $\overline{D}$ it is clear that there is a
constant $C_{1}>0$ such that
\begin{align*}
|d(p_{0},x_{n})-\log\frac{1}{a_{n}}| &  <C_{1}\\
|d(p_{0},y_{n})-\log\frac{1}{c_{n}}| &  <C_{1},
\end{align*}
for all $n$.

In view of Lemma \ref{lemma-angle} there is a constant $C_{2}>0$
such that
$$
a_{n}\leq C_{2}b_{n} \hspace{1cm} \mbox{and} \hspace{1cm} c_{n}\leq C_{2}d_{n}
$$
Therefore, for some constant $C_{3}$,
\begin{align*}
2(x_{n}|y_{n})  & =d(x_{n},p_{0})+d(y_{n},p_{0})-d(x_{n},y_{n})\\
& \geq -2C_{1}+\log\frac{1}{a_{n}}+\log\frac{1}{c_{n}}-d(x_{n},y_{n})\\
& \geq
- C_{3}+\log\frac{1}{b_{n}}+\log\frac{1}{d_{n}}-\log\frac{e_{n}+b_{n}
}{b_{n}}\frac{e_{n}+d_{n}}{d_{n}}\\
& = -C_{3}-\log(e_{n}+b_{n})(e_{n}+d_{n})\rightarrow\infty,
\end{align*}
since $e_{n}$, and at least one of $b_{n}$ and $d_{n}$ tend to
$0$.
\end{proof} \\

As an immediate application of Theorem \ref{theo-anders} we derive the
\begin{lemma} \label{lemma-accu-face}
All boundary accumulation points of a geodesic ray must belong
to one closed face. 
\end{lemma}
\begin{proof}:
Let $\gamma:[0,\infty)\rightarrow D$.
Consider two sequences $\gamma(t_i)$
and $\gamma(s_i)$ converging to the boundary. Then
\begin{displaymath}
(\gamma(t_i), \gamma(s_i))_{\gamma(0)}=\frac{1}{2}
(t_i+s_i-|t_i-s_i|)\geq\frac 12\min\{t_i,s_i\},
\end{displaymath}
which tends to infinity as $i\rightarrow\infty$.
The proposition now follows from Theorem \ref{theo-anders}.
\end{proof} \\

For two points $x,y\in D$ we denote by ${\xi}_{xy}, {\xi}_{yx}\in \partial D$ those points satisfying $x,y\in [{\xi}_{xy},{\xi}_{yx}]$
with $|{\xi}_{xy}x|<|{\xi}_{xy}y|$ and $|{\xi}_{yx},y|<|{\xi}_{yx}x|$. \\
The following lemma is well known and simply follows from the fact, that the straight line Hilbert geodesic between two points $x,y\in D$ is 
the image of the unique geodesic segment connecting the two points if and only if there do not
exist two coplanar but not collinear straight line segments $l_1,l_2\subset \partial D$ through ${\xi}_{xy}$ and ${\xi}_{yx}$ such that
${\xi}_{xy}$ and ${\xi}_{yx}$ are not boundary points of $l_1$ and $l_2$. 
\begin{lemma} \label{lemma-convergence}
Let $(D,h_D)$ be a Hilbert geometry and $x,y,z\in D$ such that $h_D(x,y)+h_D(y,z)=h_D(y,z)$ and $x\notin [y,z]$.
Then it holds $[{\xi}_{yx},{\xi}_{zy}]\subset \partial D$ and $[{\xi}_{xy},{\xi}_{yz}]\subset \partial D$.
\end{lemma}

\begin{figure}[htbp]
\psfrag{x}{$\scriptstyle{x}$} 
\psfrag{y}{$\scriptstyle{y}$}
\psfrag{z}{$\scriptstyle{z}$}
\psfrag{xizy}{$\scriptstyle{{\xi}_{zy}}$}
\psfrag{xixy}{$\scriptstyle{{\xi}_{xy}}$}
\psfrag{xiyx}{$\scriptstyle{{\xi}_{yx}}$}
\psfrag{xiyz}{$\scriptstyle{{\xi}_{yz}}$}
\psfrag{[zyyx]subsetpartialD}{$\scriptstyle{[{\xi}_{zy},{\xi}_{yx}]\subset \partial  D}$}
\psfrag{[xyyz]subsetpartialD}{$\scriptstyle{[{\xi}_{xy},{\xi}_{yz}]\subset \partial  D}$}
\psfrag{u}{$\scriptstyle{u}$}
\psfrag{v}{$\scriptstyle{v}$}
\psfrag{yn}{$\scriptstyle{y_n}$}
\psfrag{zn}{$\scriptstyle{z_n}$}
\psfrag{w}{$\scriptstyle{w}$}
\psfrag{w'}{$\scriptstyle{w'}$}
\psfrag{L(x,w)}{$\scriptstyle{L(x,w)}$}
\psfrag{xiynzn}{$\scriptstyle{{\xi}_{y_nz_n}}$}
\psfrag{xixyn}{$\scriptstyle{{\xi}_{xy_n}}$}
\psfrag{xixzn}{$\scriptstyle{{\xi}_{xz_n}}$}
\includegraphics[width=0.9\columnwidth]{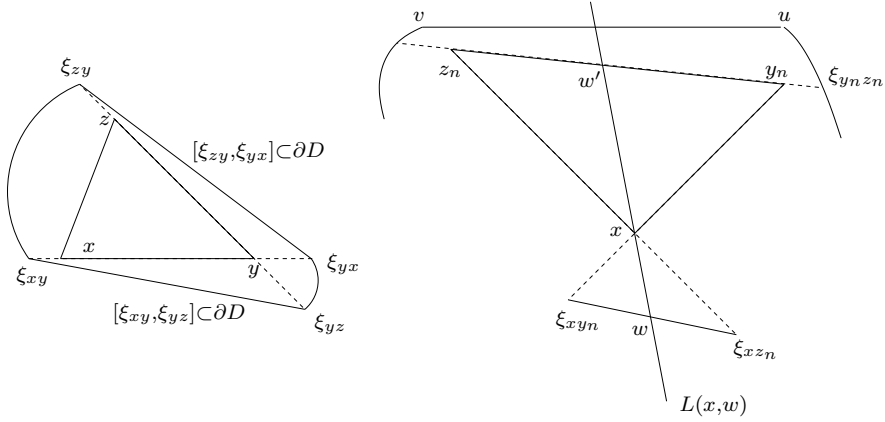}
\caption{The figure on the left hand side visualizes the situation in Lemma \ref{lemma-convergence}, while the figure on the right hand side shows the
setup of the proof of Theorem \ref{theo-convergence}.}
\end{figure}

With Proposition \ref{prop-converse} and Lemmata \ref{lemma-accu-face} and \ref{lemma-convergence} at hand, we are finally ready to provide
the \\
\begin{proof} {\bf of Theorem \ref{theo-convergence}}: 
(i) Assume to the contrary of the claim of Theorem \ref{theo-convergence} that there exist $u,v\in \partial D$, $u\neq v$, 
such that both, $u$ and $v$ are accumulation points of a geodesic ray $\gamma$ in $(D,h_D)$ with $x:=\gamma (0)$. 
From Lemma \ref{lemma-accu-face} we know that $u$ and $v$ must lie in a
common face $F$ of $\partial D$. Thus we have $[u,v]\subset \partial D$. Under this assumption we prove the following two claims, which 
cannot simultaneously hold. By that we contradict our assumption! \\
{\it Claim 1}: With the notation introduced above it holds $[{\xi}_{xu},{\xi}_{xv}]\subset \partial D$. \\
{\it Claim 2}: There exist ${\eta}_{uv}\in l:=((L(u,v)\cap \partial D)\setminus [u,v]) \cup \{ u,v\}$ with $|u{\eta}_{uv}|<|v{\eta}_{uv}|$
and ${\eta}_{vu}\in l$ with $|v{\eta}_{vu}|<|u{\eta}_{vu}|$ such that
\begin{equation} \label{eqn-eta}
[{\xi}_{xu},{\eta}_{uv}], \; [{\xi}_{xv},{\eta}_{vu}] \; \subset \; \partial D .
\end{equation}
In the following let $\{ y_n{\}}_{n\in \mathbb{N}}$,  $\{ z_n{\}}_{n\in \mathbb{N}}$ be sequences in $D$ converging to $u$ and $v$, respectively,
such that
\begin{description}
\item[(1)] $y_n,z_n\in im\{ \gamma \} \;\;\; \forall n\in \mathbb{N}$,
\item[(2)] $h_D(x,y_{n+1}) \; > \; h_D(x,z_n) \; > \; h_D(x,y_n) \;\;\; \forall n\in \mathbb{N}$,
%\item[(3)] $dist(L(u,v),y_{n+1})<dist(L(u,v)),z_n)<dist(L(u,v),y_n) \;\;\; \forall n\in \mathbb{N}$,
%where $dist(L(u,v),y_{n})$ for instance is defined through
%\begin{displaymath}
%dist\Big( L(u,v),y_{n}\Big) \; := \; \min\limits_{w\in L(u,v)} \Big\{ |y_nw|\Big\} . 
%\end{displaymath}
\end{description} 
{\bf Proof of Claim 1:} We set $E_n:=span \{ x,y_n,z_n\}$ for all $n\in \mathbb{N}$ and assume without loss of generality that
for all $n\in \mathbb{N}$ this is a nondegenerated, 2-dimensional plane. Note that the Hilbert geometry on $E_n\cap D$ isometrically
embeds into that of $D$. \\
Fix $n\in \mathbb{N}$ and suppose that $[{\xi}_{xy_n},{\xi}_{xz_n}] \not\subset \partial D$. 
Then there exists $w\in  [{\xi}_{xy_n},{\xi}_{xz_n}]$ such that all straight line segments on $L(x,w)\cap D$ are the images of the unique geodesic segments
in $(E_n\cap D,h_{E_n\cap D})$ connecting their endpoints. This clearly contradicts the fact that for $w':=L(x,w)\cap [y_n,z_n]$ the union
$[x,y_n]\cup [y_n,w']$ is the image of another geodesic segment in   $(E_n\cap D,h_{E_n\cap D})$ connecting $x$ to $w'$. Thus we have
$[{\xi}_{xy_n},{\xi}_{xz_n}]\subset \partial D$ for all $n\in \mathbb{N}$ and by 
${\xi}_{xy_n}\stackrel{n\rightarrow \infty}{\longrightarrow} {\xi}_{xu}$,
${\xi}_{xz_n}\stackrel{n\rightarrow \infty}{\longrightarrow} {\xi}_{xv}$
and the continuity of $\partial D$ also $[{\xi}_{xu},{\xi}_{xv}]\subset \partial D$. \\
{\bf Proof of Claim 2:} Let $E_n$, $n\in \mathbb{N}$, be as above and $\tilde{E}_n:=span\{ x,z_n,y_{n+1}\}$ for all $n\in \mathbb{N}$. 
As above for $E_n$ we assume without loss of generality that $\tilde{E}_n$ also is a nondegenerate, 2-dimensional plane. Now consider
the sequences $\{ {\xi}_{y_nz_n}{\}}_{n\in \mathbb{N}}$ and $\{ {\xi}_{z_ny_{n+1}}{\}}_{n\in \mathbb{N}}$. By construction we can
pass over to appropriate subsequences that converge to points in $l$. Let ${\eta}_{uv}$ and ${\eta}_{vu}$, respectively,  denote their limits.
Now it follows by the choices of the sequences $\{ y_n{\}}_{n\in \mathbb{N}}$ and $\{ z_n{\}}_{n\in \mathbb{N}}$ as well as from 
Lemma \ref{lemma-convergence} that $[{\xi}_{xy_n},{\xi}_{y_nz_n}] \subset \partial D$ and
$[{\xi}_{xy_n},{\xi}_{z_ny_{n+1}}] \subset \partial D$ for all $n\in \mathbb{N}$. Due to 
${\xi}_{xy_n}\stackrel{n\rightarrow \infty}{\longrightarrow} {\xi}_{xu}$,
${\xi}_{xz_n}\stackrel{n\rightarrow \infty}{\longrightarrow} {\xi}_{xv}$,
${\xi}_{y_nz_n}\stackrel{n\rightarrow \infty}{\longrightarrow} {\eta}_{uv}$,
${\xi}_{z_ny_{n+1}}\stackrel{n\rightarrow \infty}{\longrightarrow} {\eta}_{vu}$ and the continuity of $\partial D$ we finally obtain
the inclusions (\ref{eqn-eta}) and hence the validity of Claim 2. \\
(ii) Suppose that the geodesic line $\Gamma$ accumulates only at
$\xi\in\partial D$, so $\Gamma(t)\rightarrow\xi$ as
$t\rightarrow\pm\infty$. Then, on the one hand,
$(\Gamma(n)|\Gamma(-n))_{\Gamma (0)}\rightarrow\infty$ by 
Proposition \ref{prop-converse}, but, on the other hand, by virtue of being a geodesic,
$(\Gamma(n)|\Gamma(-n))_{\Gamma (0)}=0$. This is a contradiction. 
\end{proof} \\

Moreover, one can say that if both accumulation points of a complete geodesic 
belong to a single closed face, then both have to lie
on the boundary of the face. Such geodesic lines do exist, the
simplest example is the triangle. Here fix one point $p$ in the interior
and connect it through straight line segments to two of the vertices. The 
union of these segments indeed is the image of a geodesic in the associated
Hilbert geometry.

%%%%%%%%%%%%%%%%%%%%%%%%%%%%%%%%%%%%%%%%%%%%%%%%%%%%%%%%%%%%%%%%%%%%%%%%%%%%%%%%%%%%%%%%%%%%%
%%%%%%%%%%%%%%%%%%%%%%%%%%%%%%%%%%%%%%%%%%%%%%%%%%%%%%%%%%%%%%%%%%%%%%%%%%%%%%%%%%%%%%%%%%%%%

\end{document}